\documentclass{amsart}
\usepackage{amsfonts}
\usepackage{amsmath}
\usepackage{enumerate}
\usepackage{amsthm}
\usepackage{hyperref}
\usepackage{cite}

\begin{document}

\newtheorem{definition}{Definition}[section]
\newtheorem{theorem}[definition]{Theorem}
\newtheorem{proposition}[definition]{Proposition}
\newtheorem{remark}[definition]{Remark}
\newtheorem{lemma}[definition]{Lemma}
\newtheorem{corollary}[definition]{Corollary}
\newtheorem{example}[definition]{Example}

\numberwithin{equation}{section}

\title[Prescribed Chern scalar curvatures...]{Prescribed Chern scalar curvature on complete noncompact Hermitian manifolds with nonpositive curvatures}
\author[W. Yu]{Weike Yu}
\date{}
\thanks{The author is supported by NSFC (No. 12501075), Basic Research Program of Jiangsu (No. BK20250644), Natural Science Foundation of the Jiangsu Higher Education Institutions of China (No. 25KJB110008), the start-up research funds from Nanjing Normal University with account No. 184080H201B160, and the Laboratory of Mathematics for Nonlinear Science, Fudan University.}

\begin{abstract}
In this paper, we investigate the problem of prescribing Chern scalar curvatures on complete noncompact Hermitian manifolds with nonpositive curvatures, and establish some existence results. In particular, we obtain some sufficient conditions for the existence of constant negative Chern scalar curvature metrics in the conformal class.
\end{abstract}
\maketitle
\section{Introduction}
Suppose that $(M^n, J, h)$ is a complex manifold of complex dimension $n$ equipped with a Hermitian metric $h$ and a complex structure $J$. The fundamental form $\omega$ is defined by $\omega(\cdot, \cdot)=h(J\cdot, \cdot)$. In this paper, we shall identify the Hermitian metric $h$ with its corresponding fundamental form $\omega$. On a Hermitian manifold $(M, \omega)$, there exists a unique linear connection $\nabla^{Ch}$, called the Chern connection, preserving both the Hermitian metric $h$ and the complex structure $J$, whose torsion $T^{Ch}$ is of vanishing $(1,1)$-part everywhere. The scalar curvature with respect to $\nabla^{Ch}$, referred to as Chern scalar curvature, is given by
\begin{align}
S^{Ch}(\omega)=\text{tr}_{\omega}Ric^{(1)}(\omega)=\text{tr}_\omega \sqrt{-1}\bar{\partial}\partial  \log \omega^n,
\end{align}
where $Ric^{(1)}(\omega)$ denotes the first Chern Ricci curvature of $(M, \omega)$ and $\omega^n$ is the volume form.

A fundamental problem in differential geometry is that of studying the set of curvature functions that a given manifold can possess. The following problem is the Hermitian analogue of prescribing Riemannian scalar curvatures:
\begin{center}
\begin{minipage}{8cm}
Given a smooth real-valued function $S$ on $(M^n, \omega)$, does there exist a Hermitian metric $\tilde{\omega}$ conformal to $\omega$ such that its Chern scalar curvature $S^{Ch}(\tilde{\omega})=S$?
\end{minipage}
\end{center}
Let $\tilde{\omega}=e^{\frac{2}{n}u}\omega$ for $u\in C^\infty(M)$. Then the above problem is equivalent to solving the following partial differential equation:
\begin{align}\label{1.2}
-\Delta^{Ch}_\omega u+S^{Ch}(\omega)=Se^{\frac{2}{n}u},
\end{align}
where $\Delta^{Ch}_\omega$ denotes the Chern Laplacian associated with $\omega$. If $S$ is a constant, the above problem is called the Chern--Yamabe problem,  which was first proposed by Angella--Calamai--Spotti in \cite{[ACS]}, and they proved that it is solvable on any compact Hermitian manifold with nonpositive Gauduchon degree (see Theorem \ref{thm2.5-}). If $S$ is a general smooth function, the problem on compact (almost) Hermitian manifolds has been studied by several authors (cf. \cite{[LM], [CZ], [Ho], [Fus], [Bar], [Yu1], [WZ1], [WZ2], [LU], [LZZ]}). We remark that if $\dim_{\mathbb{C}}M=1$, the problem is known as the prescribed Gaussian curvature problem and has been extensively studied (see, for example, \cite{[AM-1], [Aub], [Ber], [BEW], [BGS], [CP], [DL], [KW1], [KW2], [KW3], [Mor]}).

In a recent work \cite{[Yu5]}, we studied the problem of prescribing Chern scalar curvature on complete Hermitian manifolds, generalizing the existence results of Aviles--McOwen \cite{[AM-1]} from Poincar\'e disks to higher-dimensional complete noncompact Hermitian manifolds whose Chern scalar curvature is bounded between two negative constants. In this paper, we further investigate this problem on complete noncompact Hermitian manifolds of nonpositive curvature. By establishing a local $C^0$-estimate for \eqref{1.2} together with the construction of lower solutions, we prove the following existence theorem, which is the Hermitian counterpart of the Aviles--McOwen existence results in \cite{[AM-2]}.

\begin{theorem}\label{theorem1.1}
Let $(M^n, \omega)$ be a complete noncompact Hermitian manifold whose second Chern Ricci curvature satisfies
\begin{align}\label{1.3}
Ric^{(2)}(X, \overline{X})\geq-C_1(1+r(x))^\alpha,
\end{align}
and whose torsion satisfies
\begin{align}\label{1.4}
\|T^{Ch}(X,Y)\|\leq C_2(1+r(x))^{\frac{\alpha}{2}},
\end{align}
for all $X, Y\in T^{1,0}_{x}M$ with $\|X\|=\|Y\|=1$. Here $r(x)$ denotes the Riemannian distance from a fixed point $x_0$ to $x$ in $(M, \omega)$, and $C_1, C_2>0, \alpha\geq 0$ are constants. Suppose further that the Chern scalar curvature of $(M, \omega)$ satisfies
\begin{align}
&S^{Ch}(\omega)(x)\leq0 , \quad\forall\ x\in M,\\
&S^{Ch}(\omega)(x)\leq -\frac{b^2}{r^l(x)}, \quad\forall\ x\in M\setminus D,
\end{align}
where $D\subset M$  is compact with $x_0$ in its interior, and $b>0$, $l\geq 0$ are constants. If either of the following two conditions holds:
\begin{enumerate}
\item $0\leq \alpha<2$ and $0\leq l<1-\frac{\alpha}{2}$;
\item $0\leq \alpha\leq 2$, $l=1-\frac{\alpha}{2}$, $b^2>4n^2\sqrt{C_3}$, where $C_3=\frac{C_1}{4n}+\frac{nC_2^2}{2}$,
\end{enumerate}
then there exists a complete Hermitian metric $\tilde{\omega}$ conformal to $\omega$ such that its Chern scalar curvature $S^{Ch}(\tilde{\omega})\equiv -1$.
\end{theorem}
\begin{remark}
In Theorem \ref{thm4.5}, we consider the problem of prescribing Chern scalar curvature on the aforementioned Hermitian manifold and establish a more general result, of which Theorem \ref{theorem1.1} is a corollary.
\end{remark}

For the case $l=0$ in Theorem \ref{theorem1.1}, however, assumptions \eqref{1.3} and \eqref{1.4} can be dropped by using a different lower-solution construction.

\begin{theorem}\label{theorem1.3}
Let $(M^n, \omega)$ be a complete noncompact Hermitian manifold whose Chern scalar curvature satisfies
\begin{align}
&S^{Ch}(\omega)(x)\leq 0, \quad\forall\ x\in M,\label{1.7+}\\ 
&S^{Ch}(\omega)(x)\leq -b^2, \quad\forall\ x\in M\setminus D,\label{1.8+}
\end{align}
where $D\subset M$ is compact and $b>0$ is a constant. Suppose $S$ is a smooth function on $M$ satisfying
\begin{align}\label{1.9}
-c^2\leq S <0.
\end{align}
for some constant $c>0$. Then there exists a complete Hermitian metric $\tilde{\omega}$ conformal to $\omega$ such that $\tilde\omega\geq C\omega$ for some constant $C>0$, and whose Chern scalar curvature satisfies $S^{Ch}(\tilde{\omega})= S$.
\end{theorem}

\begin{corollary}[\cite{[Yuan]}, \cite{[Yu5]}]
Let $(M^n, \omega)$ be a complete noncompact Hermitian manifold whose Chern scalar curvature satisfies
\begin{equation}
\begin{aligned}
&S^{Ch}(\omega)(x)\leq 0, \quad\forall\ x\in M\\
&S^{Ch}(\omega)(x)\leq -b^2, \quad\forall\ x\in M\setminus D,
\end{aligned}
\end{equation}
where $D\subset M$ is compact and $b>0$ is a constant. Then there exists a complete Hermitian metric $\tilde{\omega}$ conformal to $\omega$ such that $\tilde\omega\geq C\omega$ for some constant $C>0$, and whose Chern scalar curvature satisfies $S^{Ch}(\tilde{\omega})\equiv -1$.
\end{corollary}

This paper is organized as follows. In Section 2, we recall some basic notions and notations related to the prescribed Chern scalar curvature problem. In Section 3, we give a necessary and sufficient condition for the existence result of the prescribed Chern scalar curvature problem. In Section 4, we prove our main theorems.

\textbf{Acknowledgments.} The author would like to thank Prof. Yuxin Dong and Prof. Xi Zhang for their continued support and encouragement.

\section{Preliminaries}
In this section, we review the basic definitions and notation relevant to the prescribed Chern scalar curvature problem on Hermitian manifolds.

Let $(M^n, J, h)$ be a Hermitian manifold of complex dimension $n$, with fundamental form $\omega$ defined by $\omega(\cdot,\cdot)=h(J\cdot, \cdot)$. Let $TM^{\mathbb{C}}=TM\otimes \mathbb{C}$ be the complexification of the tangent bundle. Extending $J$ and $h$ from $TM$ to $TM^{\mathbb{C}}$ by $\mathbb{C}$-linearity, we have the decomposition:
\begin{align}
TM^{\mathbb{C}}=T^{1,0}M\oplus T^{0,1}M,
\end{align}
where $T^{1,0}M$ and $T^{0,1}M$ are the eigenspaces of the complex structure $J$ corresponding to the eigenvalues $\sqrt{-1}$ and $-\sqrt{-1}$, respectively. Furthermore, by extending the complex structure $J$ to act on differential forms, every $m$-form admits a decomposition into $(p, q)$-forms for $p, q \geq 0$ with $p + q = m$.

On a Hermitian manifold $(M^n, J, h)$, there exists a unique affine connection $\nabla^{Ch}$, called the Chern connection, which preserves both the Hermitian metric $h$ and the complex structure $J$, that is,
\begin{align}
 \nabla^{Ch}h=0,\quad \nabla^{Ch}J=0, 
\end{align}
whose torsion $T^{Ch}(X,Y)=\nabla_XY-\nabla_YX-[X,Y]$ satisfies 
\begin{align}
T^{Ch}(JX,Y)=T^{Ch}(X,JY)
\end{align}
for any $X, Y\in TM$. Throughout this paper, we will not distinguish between the Hermitian metric $h$ and its corresponding fundamental form $\omega$.

Let $(M^n, \omega)$ be a Hermitian manifold with Chern connection $\nabla^{Ch}$, and let $\{e_i\}_{i=1}^n$ be a local unitary frame of $T^{1,0}M$ with dual frame $\{\theta^i\}_{i=1}^n$. The first Chern Ricci curvature is defined by
\begin{align}
Ric^{(1)}=\sqrt{-1}Ric^{(1)}_{i\bar{j}}\theta^i\wedge\theta^{\bar{j}},\quad Ric^{(1)}_{i\bar{j}}=\sum_{k}R_{k\bar{k}i\bar{j}},
\end{align}
where $R_{i\bar{j}k\bar{l}}=h(R(e_k,e_{\bar{l}})e_i, e_{\bar{j}})$ are the components of curvature tensor $R$ of $\nabla^{Ch}$. It is well known that 
\begin{align}
Ric^{(1)}=\sqrt{-1}\bar{\partial}\partial \log{\omega^n},
\end{align}
where $\omega^n$ is the volume form of $(M^n, \omega)$. Similarly, the second Chern Ricci curvature is given by
\begin{align}
Ric^{(2)}=\sqrt{-1}Ric^{(2)}_{i\bar{j}}\theta^i\wedge\theta^{\bar{j}}, \quad Ric^{(2)}_{i\bar{j}}=\sum_{k}R_{i\bar{j}k\bar{k}}.
\end{align}
The Chern scalar curvature is then
\begin{align}
S^{Ch}(\omega)=\text{tr}_{\omega}Ric^{(1)}(\omega)=\text{tr}_{\omega}Ric^{(2)}(\omega)=\text{tr}_{\omega}\sqrt{-1}\bar{\partial}\partial \log{\omega^n}.
\end{align}

On a Hermitian manifold $(M^n, \omega)$, there is a canonical elliptic operator known as the Chern Laplacian. For any smooth function $u: M\rightarrow \mathbb{R}$, the Chern Laplacian $\Delta^{Ch}_\omega u$ is defined as
\begin{align}
\Delta^{Ch}_\omega u=-2\sqrt{-1} tr_{\omega}\overline{\partial}\partial u,
\end{align}
or, in local holomorphic coordinates $\{z^i\}$ such that $\omega\overset{\text{loc}}{=}\sqrt{-1}h_{i\bar{j}}dz^i\wedge dz^{\bar{j}}$, as
\begin{align}
\Delta^{Ch}_\omega u\overset{\text{loc}}{=}2h^{i\bar{j}}\frac{\partial^2u}{\partial z^i\partial z^{\bar{j}}},
\end{align}
where $(h^{i\bar{j}})$ is the inverse matrix of $(h_{i\bar{j}})$.
\begin{lemma}[cf. \cite{[Gau]}]\label{lemma2.1}
On a Hermitian manifold $(M^n,\omega)$, the following identity holds:
\begin{align}
-\Delta^{Ch}_\omega u=-\Delta_d u+(du,\theta)_{\omega},
\end{align}
where $\Delta_d u=-d^*du$ is the Hodge-de Rham Laplacian, and $\theta$ is the Lee form or torsion $1$-form characterized by $d\omega^{n-1}=\theta\wedge \omega^{n-1}$. Here $(\cdot, \cdot)_\omega$ stands for the induced inner product on $1$-forms. In particular, when $\omega$ is balanced (i.e., $\theta\equiv0$), we have $\Delta^{Ch}_\omega u=\Delta_d u$ for all $u\in C^\infty(M)$.
\end{lemma}
Let 
\begin{align}
\{\omega\}=\{e^{\frac{2}{n}u}\omega\ |\ u\in C^\infty(M)\}
\end{align}
denote the conformal class of the Hermitian metric $\omega$. In \cite{[Gau1]}, Gauduchon proved that
\begin{theorem}\label{theorem2.2}
Let $(M^n, \omega)$ be a compact Hermitian manifold of complex dimension $n\geq 2$. Then there exists a unique Gauduchon metric $\eta\in \{\omega\}$ (i.e., $d^*\theta=0$) with volume $1$.
\end{theorem}
In view of the above theorem, one can define an invariant $\Gamma(\{\omega\})$ of the conformal class $\{\omega\}$, called the Gauduchon degree, by
\begin{align}\label{2.6...}
\Gamma(\{\omega\})=\frac{1}{(n-1)!}\int_M c^{BC}_1(K^{-1}_M)\wedge\eta^{n-1}=\int_M S^{Ch}(\eta)d\mu_\eta,
\end{align}
where $\eta$ is the unique Gauduchon metric in $\{\omega\}$ with volume $1$, $c^{BC}_1(K^{-1}_M)$ is the first Bott-Chern class of anti-canonical line bundle $K^{-1}_M$, and $d\mu_\eta$ denotes the volume form of $\eta$. As for the relation between the sign of \eqref{2.6...} and the properties of $K^{-1}_M$ and $K_M$, we refer the reader to \cite{[Yang]}.

On a Hermitian manifold $(M^n, \omega)$, consider the conformal change 
\begin{align}
\widetilde{\omega}=e^{\frac{2}{n}u}\omega.
\end{align}
According to \cite{[Gau]}, the Chern scalar curvatures of $\widetilde{\omega}$ and $\omega$ satisfy
\begin{align}\label{2.15}
-\Delta^{Ch}_\omega u+S^{Ch}(\omega)=S^{Ch}(\widetilde{\omega})e^{\frac{2}{n}u}.
\end{align}

To resolve the prescribed Chern scalar curvature problem, it is useful to see the equation \eqref{2.15} as the Euler--Lagrange equation of some functional. Unfortunately, in general, such a functional does not always exist. Indeed, according to \cite[Prop. 5.3]{[ACS]} and \cite[Prop. 2.12]{[Fus]}, we have the following
\begin{proposition}
Equation \eqref{2.15} can be viewed as the Euler--Lagrange equation for standard $L^2$ pairing if and only if $\omega$ is balanced.
\end{proposition}

In the compact case, Angella--Calamai--Spotti \cite{[ACS]} established the following result.
\begin{theorem}\label{thm2.5-}
Let $(M^n, \omega)$ be a compact Hermitian manifold with Gauduchon degree $\Gamma(\{\omega\})\leq 0$. Then there exists a unique (up to scaling) Hermitian metric $\tilde{\omega}\in\{\omega\}$ such that its Chern scalar curvature $S^{Ch}(\tilde{\omega})=\Gamma(\{\omega\})$.
\end{theorem}
It should be noted that the case $\Gamma(\{\omega\})>0$ is more delicate and remains unsolved. 

 \section{A necessary and sufficient condition for existence}
 In this section, we prove that when the candidate curvature function is negative, equation \eqref{2.15} has a solution if and only if it admits a lower solution. We first establish a local $C^0$-estimate for solutions.
\begin{lemma}\label{lemma3.1}
Let $(M^n, \omega)$ be a complete noncompact Hermitian manifold. Suppose that $u\in C^0(M)\cap W^{1,2}_{loc}(M)$ satisfies
\begin{equation}\label{3.1}
-\Delta^{Ch}_\omega u+s= Se^{\frac{2}{n}u}
\end{equation}
in the weak sense and $u\geq u_0$ for some function $u_0\in C^0(M)$, where $s$ and $S$ are smooth functions on $M$ with $S< 0$. Then for every compact subset $X\subset M$, there exists a constant $C$ independent of $u$ such that
\begin{align}
\|u\|_{C^0(X)}\leq C.
\end{align}
\end{lemma} 
\proof Since $X$ is compact in $M$, every open cover of $X$ has a finite subcover. Hence, it suffices to prove the local estimate
\begin{align}\label{3.3..}
\|u\|_{C^0(B_\mu(x_0))}\leq C,
\end{align}
where $C>0$ is independent of $u$, and $B_{\mu}(x_0)=\{x\in M: \text{dist}_\omega(x_0, x)<\mu\}$, with $\mu>0$ so that the closed geodesic ball $\overline{B_\mu(x_0)}$ is contained within the cut locus of $x_0$ (denoted by $\text{Cut}(x_0)$). For notational convenience, we shall denote $B_\mu(x_0)$ simply by $B_\mu$.

First, we prove that there exists a constant $C>0$ independent of $u$ such that
\begin{align}\label{3.4..}
\|u\|_{W^{1,2}(B_\mu)}<C.
\end{align}
Since $u\geq u_0\in C^0(M)$, we have $u>-\|u_0\|_{C^0(\overline{B_\tau})}-1$ in $\overline{B_\tau}$ for $\tau>\mu$ and $B_\tau \cap \text{Cut}(x_0)=\emptyset$. Hence, defining $v=u+A$ with $A=\|u_0\|_{C^0(\overline{B_\tau})}+1$, we get $v>0$ on $\overline{B_\tau}$. By \cite[Lemma 3.3]{[Do]}, there exists a cutoff function $\xi\in C^\infty_0(M)$ such that
\begin{align}\label{3.12.}
\xi(x)=
\begin{cases}
1& x\in B_\mu \\
0& x\in M\setminus B_\tau
\end{cases}
,\ 0\leq\xi\leq1,\ |\nabla \xi|\leq C\xi^{\frac{1}{2}},
\end{align}
where $C>0$ is a constant and $\mu<\tau$. Taking the test function $\varphi=\xi^2 v$  in \eqref{3.1} yields
\begin{align}\label{3.6...}
\int_{B_\tau} \left(\nabla v\cdot \nabla (\xi^2 v)+\left( dv, \theta\right)_\omega \xi^2 v+ s\xi^2 v- Se^{-\frac{2}{n}A}e^{\frac{2}{n}v}\xi^2 v\right)=0.
\end{align}
Substituting \eqref{3.12.} and 
\begin{align}
\nabla v\cdot\nabla (\xi^2v)=|\nabla (\xi v)|^2-|\nabla \xi|^2v^2,\quad S(x)\leq -\epsilon_0<0\text{\ in}\ \overline{B_\tau},
\end{align}
into \eqref{3.6...}, where $\epsilon_0=-\max_{\overline{B_\tau}}S>0$ (since $S<0$ on $M$), we deduce that
\begin{equation}\label{3.8.}
\begin{aligned}
0\geq \int_{B_\tau}\left( |\nabla (\xi v)|^2-C^2\xi v^2+\left( d(\frac{1}{2}v^2),\theta\right)_\omega \xi^2  + s\xi^2 v+\epsilon_0e^{-\frac{2}{n}A}e^{\frac{2}{n}v}\xi^2 v\right).
\end{aligned}
\end{equation}
Using $|ab|\leq \delta a^2+(4\delta)^{-1}b^2\ (\forall \ \delta>0)$ and $0\leq\xi\leq1$, we have
\begin{align}
&\int_{B_\tau}\xi v^2\leq \delta_1\int_{B_\tau}\xi^2v^4+(4\delta_1)^{-1}\text{Vol}(B_{\tau}),\\
&\left|\int_{B_\tau} s\xi^2 v\right|\leq (4\delta_2)^{-1}\int_{B_\tau}s^2+\delta_2\int_{B_\tau}\xi^2v^2.
\end{align}
Besides, 
\begin{equation}\label{3.11.}
\begin{aligned}
\int_{B_\tau}\left( d(\frac{1}{2}v^2),\theta\right)_\omega \xi^2&=\frac{1}{2}\int_{B_\tau}\left( d(\xi v^2), \xi\theta\right)_\omega-\frac{1}{2}\int_{B_\tau}\xi v^2\left( d\xi, \theta\right)_\omega\\
 &=\frac{1}{2}\int_{B_\tau} \xi v^2d^*(\xi\theta)-\frac{1}{2}\int_{B_\tau}\xi v^2\left( d\xi, \theta\right)_\omega\\
 &\geq -\delta_3\int_{B_\tau}\xi v^2\\
 &\geq -\delta_3\delta_1\int_{B_\tau}\xi^2v^4-\delta_3(4\delta_1)^{-1}\text{Vol}(B_{\tau})
 \end{aligned}
 \end{equation}
where $\delta_3=\frac{1}{2}\left(\max_{\overline{B_\tau}}|d^*(\xi\theta)|+\max_{\overline{B_\tau}}|\left( d\xi, \theta\right)_\omega|\right)$. Applying the inequality $e^{2t}\geq t^3$ for any $t\in \mathbb{R}$, together with $v>0$, we obtain
\begin{equation}\label{3.12..}
\begin{aligned}
\int_{B_\tau}\epsilon_0e^{-\frac{2}{n}A}e^{\frac{2}{n}v}\xi^2 v\geq \frac{\epsilon_0}{n^3}e^{-\frac{2}{n}A}\int_{B_\tau}\xi^2 v^4.
\end{aligned}
\end{equation}
From \eqref{3.8.}-\eqref{3.12..} and the inequality $t^4>t^2-1$ for all $t\in \mathbb{R}$, it follows that
 \begin{equation}\label{3.13-}
 \begin{aligned}
 0&\geq \int_{B_\tau}|\nabla (\xi v)|^2+\left(\frac{\epsilon_0}{n^3}e^{-\frac{2}{n}A}-(C^2+\delta_3)\delta_1\right)\int_{B_\tau}\xi^2 v^4\\
& \ \ \ -\delta_2\int_{B_\tau}\xi^2v^2-(C^2+\delta_3)(4\delta_1)^{-1}\text{Vol}(B_{\tau})-(4\delta_2)^{-1}\int_{B_\tau}s^2\\
&\geq \int_{B_\tau}|\nabla (\xi v)|^2+\left(\frac{\epsilon_0}{n^3}e^{-\frac{2}{n}A}-(C^2+\delta_3)\delta_1-\delta_2\right)\int_{B_\tau}\xi^2 v^2\\
& \ \ \ -\left(\frac{\epsilon_0}{n^3}e^{-\frac{2}{n}A}-(C^2+\delta_3)\delta_1\right)\int_{B_\tau}\xi^2\\
&\ \ \ -(C^2+\delta_3)(4\delta_1)^{-1}\text{Vol}(B_{\tau})-(4\delta_2)^{-1}\int_{B_\tau}s^2\\
&=\int_{B_\theta}|\nabla v|^2+\frac{\epsilon_0}{2n^3}e^{-\frac{2}{n}A}\int_{B_\theta} v^2-n^3\epsilon_0^{-1}e^{\frac{2}{n}A}\int_{B_\tau}s^2\\
&\ \ \ -\left(\frac{3\epsilon_0}{4n^3}e^{-\frac{2}{n}A}+(C^2+\delta_3)^2n^3\epsilon_0^{-1}e^{\frac{2}{n}A} \right)\text{Vol}(B_{\tau}),
 \end{aligned}
 \end{equation}
where we have chosen $\delta_1=\frac{\epsilon_0}{4(C^2+\delta_3)n^3}e^{-\frac{2}{n}A}$, $\delta_2=\frac{\epsilon_0}{4n^3}e^{-\frac{2}{n}A}$, and used the fact that $0\leq \xi\leq 1$. Thus, from \eqref{3.13-} and $\xi\equiv 1$ in $B_{\mu}$, we have
\begin{equation}
 \begin{aligned}
 &\int_{B_\mu}|\nabla v|^2+\frac{\epsilon_0}{2n^3}e^{-\frac{2}{n}A}\int_{B_\mu} v^2\\
 &\ \ \ \leq\int_{B_\tau}|\nabla (\xi v)|^2+\frac{\epsilon_0}{2n^3}e^{-\frac{2}{n}A}\int_{B_\tau}\xi^2 v^2\\
 &\ \ \ \leq \left(\frac{3\epsilon_0}{4n^3}e^{-\frac{2}{n}A}+(C^2+\delta_3)^2n^3\epsilon_0^{-1}e^{\frac{2}{n}A} \right)\text{Vol}(B_{\tau})+n^3\epsilon_0^{-1}e^{\frac{2}{n}A}\int_{B_\tau}s^2,
  \end{aligned}
 \end{equation}
 which implies, since $u=v-A$, that
 \begin{align}
 \|u\|_{W^{1,2}(B_\mu)}\leq C'_1,
 \end{align}
for some constant $C'_1>0$ independent of $u$. 

We now establish the local $C^0$-estimate for $u$. When $\dim_{\mathbb{C}}M=n=1$, by the Trudinger--Moser inequality and \eqref{3.4..}, we have $\int_{B_\mu}e^{pu}<C'_2$ for any $p>0$, where $C'_2>0$ is a constant independent of $u$. Applying the interior elliptic $L^2$-estimates to \eqref{3.1} and using \eqref{3.4..}, we get
\begin{equation}
\begin{aligned} 
\|u\|_{W^{2,2}(B_{\frac{\mu}{2}})}&\leq C'_3\left(\|u\|_{L^2(B_\mu)}+\|Se^{2u}-s\|_{L^2(B_\mu)}\right)\\
&\leq C'_3\left(\|u\|_{L^2(B_\mu)}+\max_{\overline{B_\mu}}|S|\cdot\|e^{2u}\|_{L^2(B_\mu)}+\|s\|_{L^2(B_\mu)}\right)\\
&\leq C'_4,
\end{aligned}
\end{equation}
where $C'_3, C'_4>0$ are constants independent of $u$. Combining this with the Sobolev embedding theorem $W^{2,2}(B_{\frac{\mu}{2}})\subset C^0(\overline{B_{\frac{\mu}{2}}})$ ($\dim_{\mathbb{R}}M=2$) yields $\|u\|_{C^0(B_{\frac{\mu}{2}})}<C'_5$ for some constant $C'_5>0$ independent of $u$. 

In the case $\dim_{\mathbb{C}}M=n\geq 2$, it is enough to establish the estimate
\begin{align}\label{3.16..}
\|u\|_{C^0(B_\mu)}\leq C\max\{\|u\|_{L^2(B_\tau)},1 \},
\end{align}
with $\mu<\tau$. Once \eqref{3.16..} is proved, the desired local $C^0$-estimate follows immediately from the local $W^{1, 2}$-estimates for $u$. Take a decreasing sequence $\{r_i\}_{i=0}^\infty$ in $\mathbb{R}^+$ and choose the following cutoff functions $\xi_i\in C^\infty_0(M)$ at each step:
\begin{align}\label{3.5}
\xi_i(x)=
\begin{cases}
1& x\in B_{r_{i+1}}(x_0)\\
0& x\in M\setminus B_{r_ i}(x_0)
\end{cases}
,\ 0\leq\xi_i\leq1,\ |\nabla \xi_i|\leq \frac{C}{r_i-r_{i+1}},
\end{align}
where $r_i=\mu+\frac{\tau-\mu}{2^i}$, $\mu<r_i\leq\tau$ for $i=0,1,2, \cdots$. Since the real-valued function $x\rightarrow x|x|^{a}\ (a\geq 0)$ is differentiable with derivative $(a+1)|x|^{a}$, we are allowed to take $\varphi=\xi^2_iu|u|^{a}\in C^0(M)\cap W_0^{1,2}(M)\ (a\geq0)$ as a test function, which has compact support in $M$. Substituting this into \eqref{3.1} yields 
\begin{align}\label{3.6}
\int_M \left(\nabla u\cdot \nabla \varphi+\left( du, \theta\right)_\omega \varphi+ s\varphi- Se^{\frac{2}{n}u}\varphi\right)=0.
\end{align}
By a simple computation, we obtain
\begin{equation}\label{3.7}
\begin{aligned}
|\nabla(\xi_i u|u|^{\frac{a}{2}})|^2&=\left(\frac{a}{2}+1\right)\nabla u\cdot \nabla (\xi_i^2u|u|^{a})+|\nabla\xi_i|^2|u|^{a+2}\\
&\ \ \ \ \ \ -\frac{a}{2}\left(\frac{a}{2}+1\right)\xi_i^2| u|^a|\nabla u|^2\\ 
&\leq \left(\frac{a}{2}+1\right)\nabla u\cdot \nabla (\xi_i^2u|u|^{a})+|\nabla\xi_i|^2|u|^{a+2}.
\end{aligned}
\end{equation}
It follows from \eqref{3.5}-\eqref{3.7} that
\begin{equation}\label{3.8}
\begin{aligned}
\int_M|\nabla(\xi_i u|u|^{\frac{a}{2}})|^2&\leq \left(\frac{a}{2}+1\right)\int_M\nabla u\cdot \nabla (\xi_i^2u|u|^{a})+\int_M|\nabla\xi_i|^2|u|^{a+2}\\
&\leq \left(\frac{a}{2}+1\right)\left\{-\int_{B_{r_i}}\left( du, \theta\right)_\omega \xi_i^2u|u|^{a}-\int_{B_{r_i}} s\xi_i^2u|u|^{a}\right.\\
&\ \ \ \ \ \ \left.+\int_{B_{r_i}} Se^{\frac{2}{n}u}\xi_i^2u|u|^{a} \right\}+\frac{C^2}{(r_i-r_{i+1})^2}\int_{B_{r_i}}|u|^{a+2}.
\end{aligned}
\end{equation}
We now estimate each term in the above inequality. From \eqref{3.5}, we get
\begin{equation}\label{3.9}
\begin{aligned}
-&\int_{B_{r_i}}\left( du, \theta\right)_\omega \xi_i^2u|u|^{a}\\
&=-\frac{1}{a+2}\int_{B_{r_i}}\left( \xi_i^2d(u^2|u|^{a}), \theta\right)_\omega\\
&=-\frac{1}{a+2}\int_{B_{r_i}}\left(d(\xi_i^2u^2|u|^a), \theta\right)_\omega+\frac{1}{a+2}\int_{B_{r_i}}2\xi_iu^2|u|^a\left(d\xi_i, \theta\right)_\omega\\
&=-\frac{1}{a+2}\int_{B_{r_i}}\xi_i^2u^2|u|^ad^*\theta+\frac{1}{a+2}\int_{B_{r_i}}2\xi_iu^2|u|^a\left(d\xi_i, \theta\right)_\omega\\
&\leq \frac{1}{a+2}\max_{\overline{B_\tau}}|d^*\theta|\int_{B_{r_i}}|u|^{a+2}+\frac{1}{a+2}\int_{B_{r_i}}2|u|^{a+2}|\left(d\xi_i, \theta\right)_\omega|\\
&\leq \frac{1}{a+2}\max_{\overline{B_\tau}}|d^*\theta|\int_{B_{r_i}}|u|^{a+2}+\frac{1}{a+2}\int_{B_{r_i}}2|u|^{a+2}|d\xi_i|_\omega\cdot |\theta|_\omega\\
&\leq \frac{1}{a+2}\left(\max_{\overline{B_\tau}}|d^*\theta|+\frac{2C\max_{\overline{B_\tau}}|\theta|_\omega}{r_i-r_{i+1}}\right)\int_{B_{r_i}}|u|^{a+2},
\end{aligned}
\end{equation}
\begin{equation}\label{3.10}
\begin{aligned}
-\int_{B_{r_i}} s\xi_i^2u|u|^{a}\leq \max_{\overline{B_\tau}}|s| \int_{B_{r_i}}|u|^{a+1},
\end{aligned}
\end{equation}
\begin{equation}\label{3.11}
\begin{aligned}
\int_{B_{r_i}} Se^{\frac{2}{n}u}\xi_i^2u|u|^{a}\leq \frac{n}{2}\int_{B_{r_i}} -S|u|^a\leq \frac{n}{2}\max_{\overline{B_\tau}}(-S)\int_{B_{r_i}}|u|^a,
\end{aligned}
\end{equation}
where we have used $S<0$ and the fact $xe^x\geq -1$ for all $x\in \mathbb{R}$ in \eqref{3.11}. Substituting \eqref{3.9}-\eqref{3.11} into \eqref{3.8}, we have
\begin{equation}\label{3.12}
\begin{aligned}
\int_M|\nabla(\xi_i u|u|^{\frac{a}{2}})|^2\leq & \left(\frac{1}{2}\max_{\overline{B_\tau}}|d^*\theta|+\frac{C\max_{\overline{B_\tau}}|\theta|_\omega}{r_i-r_{i+1}}+\frac{C^2}{(r_i-r_{i+1})^2}\right)  \int_{B_{r_i}} |u|^{a+2}\\
&+\left(\frac{a}{2}+1\right)\left\{\max_{\overline{B_\tau}}|s| \int_{B_{r_i}}|u|^{a+1}+\frac{n}{2}\max_{\overline{B_\tau}}(-S)\int_{B_{r_i}}|u|^a\right\}.
\end{aligned}
\end{equation}
Using H\"older's inequality, we deduce that
\begin{equation}
\begin{aligned}\label{3.13}
&\int_{B_{r_i}}|u|^{a}\leq \frac{a}{a+2}\int_{B_{r_i}}|u|^{a+2}+\frac{2}{a+2}\text{Vol}(B_{\tau}),\\
&\int_{B_{r_i}}|u|^{a+1}\leq \frac{a+1}{a+2}\int_{B_{r_i}}|u|^{a+2}+\frac{1}{a+2}\text{Vol}(B_{\tau}).
\end{aligned}
\end{equation}
Combining \eqref{3.13} and \eqref{3.12} gives
\begin{equation}\label{3.14}
\begin{aligned}
\int_{B_{r_i}}|\nabla(\xi_i u|u|^{\frac{a}{2}})|^2\leq &C_1\left(a+1+\frac{1}{r_i-r_{i+1}}+\frac{1}{(r_i-r_{i+1})^2}\right)\int_{B_{r_i}}|u|^{a+2}\\
&+\frac{3}{2}C_1\text{Vol}(B_{\tau}),
\end{aligned}
\end{equation}
where $C_1=\max\{C^2, \sup_{\overline{B_{\tau}}}|d^*\theta|,  \max_{\overline{B_\tau}}|s|,  \frac{n}{2}\max_{\overline{B_\tau}}(-S), C\max_{\overline{B_\tau}}|\theta|_\omega\}$. By Sobolev inequality applied to $f=\xi_i u|u|^{\frac{a}{2}}\in W^{1,2}_{0}(B_\tau)$, there exists a constant $C_2=C_2(B_\tau)>0$, independent of $u$, such that
\begin{equation}\label{3.15}
\left(\int_{B_{\tau}}(\xi_i u|u|^{\frac{a}{2}})^{2\beta}\right)^{\frac{1}{\beta}}\leq C_2\int_{B_{\tau}}\left(|\nabla(\xi_i u|u|^{\frac{a}{2}})|^2+\xi_i^2 |u|^{a+2}\right),
\end{equation}
where $\beta=\frac{n}{n-1}$. By \eqref{3.5}, \eqref{3.14} and \eqref{3.15}, we have
\begin{equation}\label{3.16}
\begin{aligned}
&\left(\int_{B_{r_{i+1}}}|u|^{(a+2)\beta}\right)^{\frac{1}{\beta}}\\
&\leq C_3\left(a+2+\frac{1}{r_i-r_{i+1}}+\frac{1}{(r_i-r_{i+1})^2}\right)\int_{B_{r_i}}|u|^{a+2}+\frac{3}{2}C_3\text{Vol}(B_{\tau})\\
&\leq C_4\left(a+2+\frac{1}{r_i-r_{i+1}}+\frac{1}{(r_i-r_{i+1})^2}\right)\max\left\{\int_{B_{r_i}}|u|^{a+2}, 1\right\}
\end{aligned}
\end{equation}
where $C_3=\max\{C_1C_2, C_2\}$, $C_4=\max\{2C_3, \frac{3}{2}C_3\text{Vol}(B_{\tau})\}$. Choose $a_i$ in \eqref{3.16} such that $\beta^i=\frac{a_i}{2}+1$. Then, for each $i\geq 0$, we obtain
\begin{equation*}
\begin{aligned}
&\max\left\{\|u\|_{L^{2\beta^{i+1}}(B_{r_{i+1}})}, 1\right\}\\
& \leq C_4^{\frac{1}{2\beta^i}}\left(2\beta^i+\frac{2^{i+1}}{\tau-\mu}+\frac{4^{i+1}}{(\tau-\mu)^2}\right)^{\frac{1}{2\beta^i}}\max\left\{\|u\|_{L^{2\beta^{i}}(B_{r_{i}})}, 1\right\}\\
& = C_4^{\frac{1}{2\beta^i}}4^{\frac{i+1}{2\beta^i}}\left(\frac{1}{2}\cdot \left(\frac{\beta}{4}\right)^i+\frac{1}{2^{i+1}(\tau-\mu)}+\frac{1}{(\tau-\mu)^2}\right)^{\frac{1}{2\beta^i}}\max\left\{\|u\|_{L^{2\beta^{i}}(B_{r_{i}})}, 1\right\}\\
& \leq C_4^{\frac{1}{2\beta^i}}4^{\frac{i+1}{2\beta^i}}\left(\frac{1}{2}+\frac{1}{\tau-\mu}+\frac{1}{(\tau-\mu)^2}\right)^{\frac{1}{2\beta^i}}\max\left\{\|u\|_{L^{2\beta^{i}}(B_{r_{i}})}, 1\right\},
\end{aligned}
\end{equation*}
where $\frac{\beta}{4}=\frac{n}{4(n-1)}<1$ for all $n\geq 2$. Set $C_5=\frac{1}{2}+\frac{1}{\tau-\mu}+\frac{1}{(\tau-\mu)^2}$. Then, by iteration, it follows that
\begin{equation}
\begin{aligned}
&\max\left\{\|u\|_{L^{2\beta^{i+1}}(B_{r_{i+1}})}, 1\right\}\\
&\ \ \ \ \ \ \leq C_4^{\frac{1}{2\beta^i}}4^{\frac{i+1}{2\beta^i}}C_5^{\frac{1}{2\beta^i}}\max\left\{\|u\|_{L^{2\beta^{i}}(B_{r_{i}})}, 1\right\}\\
&\ \ \ \ \ \ \leq \cdots\\
&\ \ \ \ \ \ \leq C_4^{\frac{1}{2}\sum_{k=0}^i\frac{1}{\beta^k}}4^{\frac{1}{2}\sum_{k=0}^i\frac{k+1}{\beta^k}}C_5^{\frac{1}{2}\sum_{k=0}^i\frac{1}{\beta^k}}\max\left\{\|u\|_{L^{2}(B_{\tau})}, 1\right\}.
\end{aligned}
\end{equation}
Letting $i\rightarrow +\infty$ in the above inequality, we get
\begin{align}\label{3.31}
\|u\|_{C^0(B_\mu)}\leq \max\left\{\|u\|_{C^0(B_\mu)}, 1\right\}\leq C_6\max\{\|u\|_{L^2(B_\tau)},1 \},
\end{align}
where 
$$C_6=C_4^{\frac{1}{2}\sum_{k=0}^\infty\frac{1}{\beta^k}}4^{\frac{1}{2}\sum_{k=0}^\infty\frac{k+1}{\beta^k}}C_5^{\frac{1}{2}\sum_{k=0}^\infty\frac{1}{\beta^k}}<\infty$$
is a constant independent of $u$. By the local $W^{1,2}$-estimate \eqref{3.4..}, which holds for any $\mu$ with $B_\mu\cap \text{Cut}(x_0)=\emptyset$, we have $\|u\|_{W^{1,2}(B_\tau)}<C_7$ for some constant $C_7>0$, since $B_\tau\cap \text{Cut}(x_0)=\emptyset$. Combining this with \eqref{3.31} yields the local $C^0$-estimate \eqref{3.3..}. 
\qed

With the aid of the key local $C^0$-estimate established above, we now prove the main result of this section.

\begin{proposition}\label{prop3.2}
Let $(M^n, \omega)$ be a complete noncompact Hermitian manifold, and let $s, S$ be two smooth functions on $M$ with $S<0$. Then there exists a solution $u\in C^\infty(M)$ of 
\begin{align}
-\Delta^{Ch}_\omega u+s=Se^{\frac{2}{n}u} \quad\text{in $M$}
\end{align}
if and only if there exists a function $u_-\in C^0(M)\cap W^{1,2}_{loc}(M)$ such that
\begin{align}\label{3.4}
-\Delta^{Ch}_\omega u_-+s\leq Se^{\frac{2}{n}u_-} \quad\text{in $M$}
\end{align}
holds in the weak sense. Moreover, $u\geq u_-$ in $M$.
\end{proposition}
\proof  It is clear that the ``only if" direction is trivial. Now suppose there exists a desired lower solution $u_-\in C^0(M)\cap W^{1,2}_{loc}(M)$ satisfying \eqref{3.4}. Let $M=\cup_{k\geq 1} \Omega_k$, where $\Omega_k$ is bounded domain in $M$ with smooth boundary satisfying $\Omega_k\subset\subset \Omega_{k+1}$. For each $\Omega_k$, since $s, S\in C^\infty(M)$ with $S<0$ on $M$, there exist two positive constants $C_k, C'_k$ depending on $\Omega_k$ such that $s>-C_k$ and $S\leq -C'_k<0$ in $\Omega_k$. We then have an upper solution $u_+=b_k$ in $\Omega_k$, where $b_k>\max\{\frac{n}{2}\log \frac{C_k}{C'_k}, \max_{\overline{\Omega}_k}u_-\}$ is a constant. Indeed,
\begin{align}
-\Delta^{Ch}_{\omega} u_++s-Se^{\frac{2}{n}u_+}=s-Se^{\frac{2}{n}b_k}>-C_k+C'_ke^{\frac{2}{n}b_k}>0\ \text{in}\ \Omega_k.
\end{align}
By the Monotone Iteration Scheme (see the proof of \cite[Prop. 4.1]{[Yu5]} for details), there is a solution $u_k\in C^\infty(\Omega_k)$ of 
\begin{equation}\label{4.68}
\begin{cases}
-\Delta^{Ch}_{\omega} u_k+s=Se^{\frac{2}{n}u_k}\quad \text{in}\ \Omega_k,\\
u_-\leq u_k\leq u_+=b_k\quad \text{in}\ \Omega_k.
\end{cases}
\end{equation}

Consider the sequence $\{u_k\}_{k\geq 4}$ on $\overline{\Omega}_3$. By Lemma \ref{lemma3.1}, we have $\|u_k\|_{C^0(\overline{\Omega}_3)}\leq C$ for some constant $C>0$ independent of $k$. By the interior elliptic $L^p$-estimates, there is a constant $C>0$ such that for any $p>1$,
\begin{equation*}
\begin{aligned}
\|u_k\|_{W^{2,p}(\Omega_2)}&\leq C\left(\|u_k\|_{L^p(\Omega_3)}+\|Se^{\frac{2}{n}u_k}-s\|_{L^p(\Omega_3)}\right)\\
&\leq  C\left(\|u_k\|_{C^0(\overline{\Omega}_3)}|\Omega_3|^{\frac{1}{p}}+\|S\|_{C^0(\overline{\Omega}_3)}e^{\frac{2}{n}\|u_k\|_{C^0(\overline{\Omega}_3)}}|\Omega_3|^{\frac{1}{p}}+\|s\|_{L^p(\Omega_3)}\right)\\
&\leq C'
\end{aligned}
\end{equation*}
where $C'$ is a constant independent of $k$. Taking $p=2n+1$ and applying the Sobolev embedding theorem $W^{2, 2n+1}(\Omega_2)\subset C^1(\Omega_2)$ yields $\|u_k\|_{C^1(\Omega_2)}\leq C$. By the Schauder estimates, we have $\|u_k\|_{C^{2, \alpha}(\Omega_1)}\leq C$ for any $\alpha\in (0,1)$. The compact embedding $C^{2, \alpha}(\Omega_1)\subset\subset C^{2}(\Omega_1)$ then implies the existence of a subsequence $\{u_{k1}\}\subset \{u_k\}$ such that $u_{k1}\rightarrow u_1$ in $C^2(\Omega_1)$ as $k\rightarrow +\infty$, where $u_1$ is a solution of
\begin{align}
-\Delta^{Ch}_{\omega} u+s=Se^{\frac{2}{n}u}\quad \text{in}\ \Omega_1.
\end{align}
Repeating the above procedure with the sequence $\{u_{k1}\}_{k\geq 5}$ on $\overline{\Omega}_4$, there exists a subsequence $\{u_{k2}\}\subset \{u_{k1}\}$ such that  $u_{k2}\rightarrow u_2$ in $C^2(\Omega_2)$ as $k\rightarrow +\infty$, and $u_2$ solves
\begin{equation}
\begin{cases}
-\Delta^{Ch}_{\omega} u_2+s=Se^{\frac{2}{n}u_2}\quad \text{in}\ \Omega_2,\\
u_2|_{\Omega_1}=u_1.
\end{cases}
\end{equation}
Inductively, we obtain a subsequence $\{u_{kj}\}$ which converges to $u_j$ in $C^2(\Omega_j)$ as $k\rightarrow +\infty$. Here $u_j$ is a solution of 
\begin{equation}
\begin{cases}
-\Delta^{Ch}_{\omega} u_j+s=Se^{\frac{2}{n}u_j}\quad \text{in}\ \Omega_j,\\
u_j|_{\Omega_{j-1}}=u_{j-1}
\end{cases}
\end{equation}
for any $j\geq 2$. Define 
\begin{align}
u(x)=\lim_{k\rightarrow +\infty}u_{kk}(x)\quad \ x\in M.
\end{align}
Then $u\in C^2(M)$ and it solves
\begin{equation}
-\Delta^{Ch}_{\omega} u+s=Se^{\frac{2}{n}u}\quad \text{in}\ M.
\end{equation}
By the Schauder estimates, we obtain $u\in C^\infty(M)$. Moreover, from \eqref{4.68}, it follows that $u\geq u_-$ on $M$. This completes the proof of the proposition.
\qed

\section{Construction of lower solutions and proof of the main theorems}
In this section, we prove the main theorems of this paper by constructing lower solutions and applying Proposition \ref{prop3.2}. To this end, we first recall the following Chern Laplacian comparison theorem.
\begin{lemma}[\cite{[Yu5]}]\label{lemma4.1}
Let $(M^n, \omega)$ be a complete noncompact Hermitian manifold whose second Chern Ricci curvature satisfies
\begin{align}\label{4.1}
Ric^{(2)}(X, \overline{X})\geq-C_1(1+r(x))^\alpha,
\end{align}
and whose torsion satisfies
\begin{align}\label{4.2}
\|T^{Ch}(X,Y)\|\leq C_2(1+r(x))^\beta,
\end{align}
for all $X, Y\in T^{1,0}_{x}M$ with $\|X\|=\|Y\|=1$. Here $r(x)$ denotes the Riemannian distance from a fixed point $x_0$ to $x$ in $(M, \omega)$, and $C_1, C_2>0$, $\alpha, \beta\geq 0$ are constants. Define 
\begin{align}
&F(r)=\frac{C_1}{4n}(1+r)^\alpha+\frac{nC_2^2}{2}(1+r)^{2\beta},\label{4.3}\\
&h(t)=\frac{1}{\sqrt{F(0)}}\left(e^{\int^t_0F^{\frac{1}{2}}(s)ds}-1\right),\quad t\geq 0.\label{4.4}
\end{align}
Then we have
\begin{align}\label{4.5}
\Delta^{Ch}_\omega r(x)\leq  4n \frac{h'(r(x))}{h(r(x))}
\end{align}
for any $x\in M\setminus \left(\{x_0\}\cup \text{Cut}(x_0)\right)$, where $\text{Cut}(x_0)$ denotes the cut locus of $x_0$.
\end{lemma}
\begin{remark}
In \cite{[Yu5]}, we proved the above Chern Laplacian comparison theorem for $\alpha, \beta>0$. It is clear that our proof also works for $\alpha, \beta\geq 0$.
\end{remark}

If we take $\beta=\frac{\alpha}{2}\geq 0$ in the above lemma, then from \eqref{4.3} we obtain
\begin{align}
F(r)=C_3(1+r)^\alpha,
\end{align}
where $C_3=\frac{C_1}{4n}+\frac{nC_2^2}{2}$. A straightforward computation using \eqref{4.4} and \eqref{4.5} then yields

\begin{corollary}\label{coro4.3}
Let $(M^n, \omega)$ be a complete noncompact Hermitian manifold whose second Chern Ricci curvature satisfies
\begin{align}
Ric^{(2)}(X, \overline{X})\geq-C_1(1+r(x))^\alpha,
\end{align}
and whose torsion satisfies
\begin{align}
\|T^{Ch}(X,Y)\|\leq C_2(1+r(x))^{\frac{\alpha}{2}},
\end{align}
for all $X, Y\in T^{1,0}_{x}M$ with $\|X\|=\|Y\|=1$. Here $r(x)$ denotes the Riemannian distance from a fixed point $x_0$ to $x$ in $(M, \omega)$, and $C_1, C_2>0$, $\alpha\geq 0$ are constants. Then we have
\begin{align}\label{4.9+}
\Delta^{Ch}_\omega r(x)\leq  \frac{2n(\alpha+2)}{r(x)}+4n\sqrt{C_3}(1+r(x))^{\frac{\alpha}{2}}
\end{align}
for any $x\in M\setminus \left(\{x_0\}\cup \text{Cut}(x_0)\right)$, where $\text{Cut}(x_0)$ denotes the cut locus of $x_0$.
\end{corollary}

We now prove the main theorem of this paper. Indeed, using the above corollary, we construct a lower solution $u_-\in C^0(M)\cap W^{1,2}_{loc}(M)$ of 
\begin{align}\label{4.10...}
-\Delta^{Ch}_\omega u+s=Se^{\frac{2}{n}u} 
\end{align}
in the weak sense, under certain assumptions on the curvature. By Proposition \ref{prop3.2}, this implies the existence of a solution $u\in C^\infty(M)$ to \eqref{4.10...}.

\begin{proposition}\label{thm4.4}
Let $(M^n, \omega)$ be a complete noncompact Hermitian manifold whose second Chern Ricci curvature satisfies
\begin{align}
Ric^{(2)}(X, \overline{X})\geq-C_1(1+r(x))^\alpha,
\end{align}
and whose torsion satisfies
\begin{align}
\|T^{Ch}(X,Y)\|\leq C_2(1+r(x))^{\frac{\alpha}{2}},
\end{align}
for all $X, Y\in T^{1,0}_{x}M$ with $\|X\|=\|Y\|=1$, where $r(x)$ is the Riemannian distance from a fixed point $x_0$ to $x$ in $(M, \omega)$, $C_1, C_2>0$ and $\alpha\in [0, 2]$. Suppose that the Chern scalar curvature of $(M, \omega)$ satisfies
\begin{align}
&S^{Ch}(\omega)(x)\leq0 , \quad\forall\ x\in M,\label{4.13--}\\
&S^{Ch}(\omega)(x)\leq -\frac{b^2}{r^l(x)}, \quad\forall\ x\in M\setminus D,\label{4.14-}
\end{align}
where $D\subset M$ is compact with $x_0$ in its interior, and $b>0,\ l\geq0$ are constants. Suppose that $S$ is a smooth function on $M$ satisfying 
\begin{align}\label{4.15...}
-c^2(1+r^k(x))\leq S(x) <0. \quad \forall\ x\in M,
\end{align}
for some constants $c>0,\ k\geq 0$. If one of the following four conditions holds:
\begin{enumerate}
\item $0\leq \alpha<2$, $0\leq l<1-\frac{\alpha}{2}$, $0\leq k<2-l$;
\item $0\leq \alpha<2$, $0\leq l<1-\frac{\alpha}{2}$, $k=2-l$, $c<b$;
\item $0\leq \alpha\leq2$, $l=1-\frac{\alpha}{2}$, $k=1+\frac{\alpha}{2}$, $c^2+8n^2\sqrt{C_3}<b^2$;
\item $0\leq \alpha\leq2$, $l=1-\frac{\alpha}{2}$, $0\leq k<1+\frac{\alpha}{2}$, $8n^2\sqrt{C_3}<b^2$,
\end{enumerate}
where $C_3=\frac{C_1}{4n}+\frac{nC_2^2}{2}$, then there exists a complete Hermitian metric $\tilde{\omega}$ conformal to $\omega$ such that its Chern scalar curvature $S^{Ch}(\tilde{\omega})= S$.
\end{proposition}
\proof 
Let $D\subset D_1\subset\subset D_2\subset\subset M$, where $D_1$ and $D_2$ are relatively compact open subsets of $M$ with smooth boundaries. Let $\phi_0\in C^\infty(\overline{D_2})$ be the solution of 
\begin{align}\label{4.16+}
\begin{cases}
&\Delta^{Ch}_\omega \phi_0=\delta_1 \quad\text{in}\ D_2,\\
& \phi_0|_{\partial D_2}=0,
\end{cases}
\end{align}
for some constant $\delta_1>0$. Note that by applying the maximum principle to \eqref{4.16+}, we obtain $\phi_0<0$ in $D_2$. Extend $\phi_0|_{D_1}$ smoothly to a compactly supported function on $M$, and denote the resulting function by $\phi_1$. Then $\phi_1\in C^\infty_0(M)$ satisfies $\Delta^{Ch}_\omega \phi_1=\delta_1$ on $D$ and $\Delta^{Ch}_\omega \phi_1\geq-C$ in $M$ for some constant $C>0$. By \eqref{4.14-}, we have
\begin{align}
S^{Ch}(\omega)(x)\leq -\frac{b^2}{\max_{\text{supp}\phi_1}r^l}=:-\tilde{b}^2<0 , \quad\forall\ x\in \text{supp} \phi_1\setminus D.
\end{align}
Set $\phi=\frac{\tilde{b}^2}{2C}\phi_1\in  C^\infty_0(M)$. Then $\Delta^{Ch}_\omega \phi=\delta$ on $D$ and $\Delta^{Ch}_\omega \phi\geq -\frac{\tilde{b}^2}{2}$ in $M$, where $\delta=\frac{\tilde{b}^2}{2C}\delta_1>0$. Let $\omega_1=e^{\frac{2}{n}\phi}\omega$. Then we have
\begin{align}
&S^{Ch}(\omega_1)=e^{-\frac{2}{n}\phi}\left(-\Delta^{Ch}_\omega\phi+S^{Ch}(\omega)\right)\leq -e^{-\frac{2}{n}\|\phi\|_{C^0}}\delta<0\quad \text{in}\ D,\label{4.18...}\\
&S^{Ch}(\omega_1)\leq e^{-\frac{2}{n}\phi}\left(\frac{\tilde{b}^2}{2}-\tilde{b}^2\right)\leq -e^{-\frac{2}{n}\|\phi\|_{C^0}}\frac{\tilde{b}^2}{2}<0\quad \text{in\ supp\ }\phi \setminus D,\\
&S^{Ch}(\omega_1)=S^{Ch}(\omega)\leq -\frac{b^2}{r^l(x)} \quad \text{in\ } M\setminus \text{supp\ }\phi.\label{4.20...}
\end{align}
If $\tilde{\omega}=e^{\frac{2}{n}u}\omega_1$ for some $u\in C^\infty(M)$, then $u$ satisfies
\begin{align}
-\Delta^{Ch}_{\omega_1} u+S^{Ch}(\omega_1)=S^{Ch}(\tilde{\omega})e^{\frac{2}{n}u}.
\end{align}
By Proposition \ref{prop3.2}, it suffices to prove the existence of a lower solution $u_-\in C^0(M)\cap W^{1,2}_{loc}(M)$ such that
\begin{align}
-\Delta^{Ch}_{\omega_1} u_-+s\leq Se^{\frac{2}{n}u_-} \quad\text{in $M$}
\end{align}
holds in the weak sense, where $s=S^{Ch}(\omega_1)$. Consider the function
\begin{align}\label{4.23+}
u_-=-\frac{n}{2}\log{(r^2+a)}\in C^0(M)\cap W^{1,2}_{loc}(M),
\end{align}
where $a>0$ is a constant to be chosen later. A direct computation gives, in $M\setminus (\{x_0\}\cup\text{Cut}(x_0))$,
\begin{equation}\label{4.24+}
\begin{aligned}
&-\Delta^{Ch}_{\omega_1} u_-+s-Se^{\frac{2}{n}u_-}\\
&=\frac{n(a-r^2)|\nabla r|^2_{\omega_1}+n(r^2+a)r\Delta^{Ch}_{\omega_1} r}{(r^2+a)^2}+s-\frac{S}{r^2+a}\\
&=\frac{n(a-r^2)e^{-\frac{2}{n}\phi}|\nabla r|^2_{\omega}+n(r^2+a)re^{-\frac{2}{n}\phi}\Delta^{Ch}_{\omega} r}{(r^2+a)^2}+s-\frac{S}{r^2+a}.
\end{aligned}
\end{equation}
Using Corollary \ref{coro4.3}, \eqref{4.15...} and \eqref{4.20...}, we have, in $M\setminus (\text{supp}\phi\cup \text{Cut}(x_0))$,
\begin{equation}\label{4.25+}
\begin{aligned}
&-\Delta^{Ch}_{\omega_1} u_-+s-Se^{\frac{2}{n}u_-}\\
=&\frac{n(a-r^2)|\nabla r|^2_{\omega}+n(r^2+a)r\Delta^{Ch}_{\omega} r}{(r^2+a)^2}+s-\frac{S}{r^2+a}\\
\leq &\frac{n(a-r^2)+n(r^2+a)\left(2n(\alpha+2)+4n\sqrt{C_3}r(1+r)^{\frac{\alpha}{2}}\right)}{(r^2+a)^2}\\
&\ \ \  -b^2r^{-l}+\frac{c^2(1+r^{k})}{r^2+a}\\
\leq& (r^2+a)^{-2}\left \{n\left(a-r^2\right)+n\left(r^2+a\right)\left(2n(\alpha+2)+4n\sqrt{C_3}r(1+r)^{\frac{\alpha}{2}}\right)\right.\\
&\ \ \ \left.-b^2r^{-l}\left(r^2+a\right)^2+c^2\left(1+r^{k}\right)\left(r^2+a\right)\right\}.
\end{aligned}
\end{equation}
Thus, $-\Delta^{Ch}_{\omega_1} u_-+s-Se^{\frac{2}{n}u_-}\leq 0$ holds on $M\setminus (\text{supp}\phi\cup \text{Cut}(x_0))$, provided that the right-hand side of \eqref{4.25+} is nonpositive, i.e.,
\begin{equation}\label{4.18}
\begin{aligned}
0&\geq n(a-r^2)+n(r^2+a)(2n(\alpha+2)+4n\sqrt{C_3}r(1+r)^{\frac{\alpha}{2}})\\
&\ \ \ \ -b^2r^{-l}(r^2+a)^2+c^2(1+r^{k})(r^2+a)\\
&=-b^2r^{4-l}-2ab^2r^{2-l}-a^2b^2r^{-l}+c^2r^{2+k}+ac^2r^k\\
&\ \ \ \ +4n^2\sqrt{C_3}(1+r)^{\frac{\alpha}{2}}(r^3+ar)+\left(2n^2(\alpha+2)+c^2-n\right)r^2\\
&\ \ \ \ +\left(na+2n^2(\alpha+2)a+ac^2 \right).
\end{aligned}
\end{equation}
If one of the four conditions $(1)-(4)$ holds, it follows from \eqref{4.18} that there exists a constant $r_0>0$ such that $\text{supp}\phi\subset B_{r_0}(x_0)=\{x\in M:\ r(x)<r_0\}$ and
\begin{align}
-\Delta^{Ch}_{\omega_1} u_-+s-Se^{\frac{2}{n}u_-}\leq 0
\end{align}
holds on $ M\setminus (\overline{B_{r_0}(x_0)}\cup\text{Cut}(x_0))$. Indeed, if condition (1) holds, i.e., 
\begin{align}
0\leq \alpha<2,\quad 0\leq l<1-\frac{\alpha}{2},\quad 0\leq k<2-l,
\end{align}
then the right-hand side of \eqref{4.18} tends to $-\infty$ as $r(x)\rightarrow +\infty$. Therefore, there exists a constant $r_0>0$ such that \eqref{4.18} holds for any $x\in M$ with $r(x)>r_0$. The remaining three cases are proved in a similar manner, so we omit it.

By \eqref{4.18...}-\eqref{4.20...}, we have $s=S^{Ch}(\omega_1)<0$ on $M$. Then, for any $x\in \overline{B_{r_0}(x_0)}\setminus (\text{Cut}(x_0)\cup\{x_0\})$, it follows from \eqref{4.24+}, \eqref{4.9+} \eqref{4.15...} that
\begin{equation}
\begin{aligned}
&-\Delta^{Ch}_{\omega_1} u_-(x)+s(x)-S(x)e^{\frac{2}{n}u_-(x)}\\
\leq &e^{\frac{2}{n}\|\phi\|_{C^0(M)}}\frac{n(a-r^2)+n(r^2+a)r\Delta^{Ch}_{\omega} r}{(r^2+a)^2}+\max_{\overline{B_{r_0}(x_0)}}s+\frac{c^2(1+r^{k})}{r^2+a}\\
\leq & e^{\frac{2}{n}\|\phi\|_{C^0(M)}}\frac{na+n(r^2_0+a)\left(2n(\alpha+2)+4n\sqrt{C_3}r_0(1+r_0)^{\frac{\alpha}{2}}\right)}{a^2}\\
&\ \ \ +\max_{\overline{B_{r_0}(x_0)}}s+\frac{c^2(1+r_0^{k})}{a}\\
\rightarrow & \max_{\overline{B_{r_0}(x_0)}}s<0, \quad \text{as\ } a\rightarrow +\infty. 
\end{aligned}
\end{equation}
Hence, we can choose a constant $a>0$ sufficiently large such that
\begin{align}
-\Delta^{Ch}_{\omega_1} u_-(x)+s(x)-S(x)e^{\frac{2}{n}u_-(x)}\leq 0
\end{align}
holds for any $x\in \overline{B_{r_0}(x_0)}\setminus (\text{Cut}(x_0)\cup\{x_0\})$. 

Consequently, we have proved that $u_-\in C^0(M)\cap W^{1,2}_{loc}(M)\cap C^\infty(M\setminus (\text{Cut}(x_0)))$ satisfies
\begin{align}\label{4.33....}
-\Delta^{Ch}_{\omega_1} u_-+s-Se^{\frac{2}{n}u_-}\leq 0 \quad \text{in\ } M\setminus (\text{Cut}(x_0)\cup\{x_0\}), 
\end{align}
where the measure $|\text{Cut}(x_0)|=0$. We now prove that
\begin{align}
-\Delta^{Ch}_\omega u_-+s\leq Se^{\frac{2}{n}u_-} \quad\text{in $M$}
\end{align}
holds in the weak sense. Let $E_0$ be the maximal domain of definition of normal coordinates at $x_0$; thus $\exp_{x_0}|_{E_0}$ is a diffeomorphism onto its image and we have $\text{Cut}(x_0)=\partial (\exp_{x_0}(E_0))$. Since $E_0$ is a star-shaped domain, we can exhaust $E_0$ by a family $\{E_0^n\}$ of relatively compact, star-shaped domains with smooth boundary. We set $\Omega^n=\exp_{x_0}(E^n_0)$, $\Omega=M\setminus \text{Cut}(x_0)$ so that
\begin{align}
\overline{\Omega^n}\subset \Omega^{n+1}, \quad \bigcup_n\Omega^n=\Omega.
\end{align}
Since each $E^n_0$ is star-shaped, we have
\begin{align}\label{4.36+}
\frac{\partial r}{\partial \nu_n}>0 \quad \text{on\ } \partial\Omega^n,
\end{align}
where $\nu_n$ is the outward unit normal to $\partial \Omega^n$. Let $0\leq \varphi\in C^n_0(M)$ and let $\delta>0$ be a small constant. By \eqref{4.33....} and the second Green formula, we obtain
\begin{equation}\label{4.37....}
\begin{aligned}
0&\geq \int_{\Omega^n\setminus B_\delta(x_0)} \left(-\Delta_{d} u_-+\langle du_-, \theta\rangle+s-Se^{\frac{2}{n}u_-}\right)\varphi\\
&=\int_{\Omega^n\setminus B_\delta(x_0)} -u_-\Delta_d\varphi-\int_{\partial \Omega^n\cup \partial B_\delta(x_0)}\left(\varphi\frac{\partial u_-}{\partial \nu_n}-u_-\frac{\partial \varphi}{\partial \nu_n}\right)\\
&\ \ \ +\int_{\Omega^n\setminus B_\delta(x_0)}  \left(\langle du_-, \theta\rangle+s-Se^{\frac{2}{n}u_-}\right)\varphi.
\end{aligned}
\end{equation}
According to \eqref{4.23+} and \eqref{4.36+}, we obtain
\begin{align}
\frac{\partial u_-}{\partial \nu_n}=-\frac{nr}{r^2+a}\frac{\partial r}{\partial \nu_n}<0\quad \text{on\ }\partial \Omega^n,
\end{align}
which implies
\begin{align}\label{4.39....}
\int_{ \partial \Omega^n}\varphi\frac{\partial u_-}{\partial \nu_n}\leq 0.
\end{align}
Submitting \eqref{4.39....} into \eqref{4.37....} yields
\begin{equation}\label{4.40+}
\begin{aligned}
\int_{\Omega^n}  \left( -u_-\Delta_d\varphi+\langle du_-, \theta\rangle\varphi+s\varphi-Se^{\frac{2}{n}u_-}\varphi\right)\leq I_\delta+\epsilon_n,
\end{aligned}
\end{equation}
where 
\begin{align*}
&I_\delta=\int_{B_\delta}  \left( -u_-\Delta_d\varphi+\langle du_-, \theta\rangle\varphi+s\varphi-Se^{\frac{2}{n}u_-}\varphi\right)+\int_{ \partial B_\delta}\left(\varphi\frac{\partial u_-}{\partial r}-u_-\frac{\partial \varphi}{\partial r}\right),\\
&\epsilon_n=-\int_{ \partial \Omega^n}u_-\frac{\partial \varphi}{\partial \nu_n}.
\end{align*}
Clearly, $I_\delta\rightarrow 0$ as $\delta\rightarrow 0+$. On the other hand, since $\varphi\in C^\infty_0(M)$ and $|\text{Cut}(x_0)|=0$, applying the divergence and Lebesgue theorems gives
\begin{align}
\epsilon_n=\int_{\Omega^n}\text{div}(u_-\nabla \varphi)\rightarrow \int_{\Omega}\text{div}(u_-\nabla \varphi)=\int_{M}\text{div}(u_-\nabla \varphi)=0,
\end{align}
as $n\rightarrow +\infty$. Now, letting $\delta\rightarrow 0+$ and $n\rightarrow +\infty$ in \eqref{4.40+} and applying integration by parts, we obtain
\begin{equation}
\begin{aligned}
0&\geq \int_{\Omega}  \left( -u_-\Delta_d\varphi+\langle du_-, \theta\rangle\varphi+s\varphi-Se^{\frac{2}{n}u_-}\varphi\right)\\
&= \int_{M}  \left( -u_-\Delta_d\varphi+\langle du_-, \theta\rangle\varphi+s\varphi-Se^{\frac{2}{n}u_-}\varphi\right)\\
&= \int_{M}  \left( \nabla u_-\cdot\nabla\varphi+\langle du_-, \theta\rangle\varphi+s\varphi-Se^{\frac{2}{n}u_-}\varphi\right).
\end{aligned}
\end{equation}
Therefore, $u_-\in C^0(M)\cap W^{1,2}_{loc}(M)$ satisfies
\begin{align}
-\Delta^{Ch}_\omega u_-+s\leq Se^{\frac{2}{n}u_-} \quad\text{in $M$}
\end{align}
in the weak sense. By Proposition \ref{prop3.2}, there exists a solution $u\in C^\infty(M)$ of 
\begin{align}
-\Delta^{Ch}_\omega u+s=Se^{\frac{2}{n}u} \quad\text{in $M$}
\end{align}
with $u\geq u_-$. Moreover, for any $x\in M$ with $r(x)>\sqrt{a}$, we have, at $x$,
\begin{align}
\tilde{\omega}=e^{\frac{2}{n}u}\omega_1\geq e^{\frac{2}{n}u_-}\omega_1=\frac{1}{r^2(x)+a}e^{\frac{2}{n}\phi}\omega\geq \frac{e^{-\frac{2}{n}\|\phi\|_{C^0(M)}}}{2r^2(x)}\omega.
\end{align}
Thus, by \cite[Lemma 5.2]{[AM-2]}, the metric $\tilde{\omega}$ is complete as well. This completes the proof of the proposition.
\qed

The above proposition yields the following existence results, which constitute one of the main theorems of this paper.

\begin{theorem}\label{thm4.5}
Let $(M^n, \omega)$ be a complete noncompact Hermitian manifold whose second Chern Ricci curvature satisfies
\begin{align}
Ric^{(2)}(X, \overline{X})\geq-C_1(1+r(x))^\alpha,
\end{align}
and whose torsion satisfies
\begin{align}
\|T^{Ch}(X,Y)\|\leq C_2(1+r(x))^{\frac{\alpha}{2}},
\end{align}
for all $X, Y\in T^{1,0}_{x}M$ with $\|X\|=\|Y\|=1$, where $r(x)$ denotes the Riemannian distance from a fixed point $x_0$ to $x$ in $(M, \omega)$, and $C_1, C_2>0, \alpha\geq 0$ are constants. Suppose that the Chern scalar curvature of $(M, \omega)$ satisfies
\begin{align}
&S^{Ch}(\omega)(x)\leq0 , \quad\forall\ x\in M,\\
&S^{Ch}(\omega)(x)\leq -\frac{b^2}{r^l(x)}, \quad\forall\ x\in M\setminus D,\label{4.14--}
\end{align}
where $D\subset M$ is compact with $x_0$ in its interior, and $b>0,\ l\geq0$ are constants. Suppose that $\tilde{S}$ is a smooth function on $M$ satisfying 
\begin{align}\label{4.33.}
-\tilde{c}^2(1+r^k(x))\leq \tilde{S}(x) <0, \quad \forall\ x\in M,
\end{align}
for some constants $\tilde{c}>0,\ k\geq 0$. If one of the following two conditions holds:
\begin{enumerate}
\item $0\leq \alpha<2$, $0\leq l<1-\frac{\alpha}{2}$, $0\leq k\leq 2-l$;
\item $0\leq \alpha\leq 2$, $l=1-\frac{\alpha}{2}$, $0\leq k\leq 1+\frac{\alpha}{2}$, $b^2>4n^2\sqrt{C_3}$,
\end{enumerate}
where $C_3=\frac{C_1}{4n}+\frac{nC_2^2}{2}$, then there exists a complete Hermitian metric $\tilde{\omega}$ conformal to $\omega$ such that its Chern scalar curvature $S^{Ch}(\tilde{\omega})= \tilde{S}$.
\end{theorem}
\proof If 
\begin{align}
0\leq \alpha<2, \quad 0\leq l<1-\frac{\alpha}{2},\quad 0\leq k< 2-l,
\end{align}
 or 
 \begin{align}
 0\leq \alpha\leq 2, \quad l=1-\frac{\alpha}{2},\quad 0\leq k< 1+\frac{\alpha}{2},\quad b^2>4n^2\sqrt{C_3},
 \end{align}
 then the result of this theorem follows from Proposition \ref{thm4.4}. If
 \begin{align}
 0\leq \alpha<2, \quad 0\leq l<1-\frac{\alpha}{2},\quad k= 2-l,
 \end{align}
 then we set 
\begin{align}\label{4.36-}
S=e^{\frac{2}{n}a}\tilde{S}, \quad \quad c^2=e^{\frac{2}{n}a}\tilde{c}^2,
\end{align}
where $a<\frac{n}{2}\log{\frac{b^2}{\tilde{c}^2}}$ is a constant. From \eqref{4.33.} and \eqref{4.36-}, we have
\begin{align}
-c^2(1+r^k)\leq S<0, \quad c^2<b^2.
\end{align}
Thus, by Proposition \ref{thm4.4} (2), there exists a complete Hermitian metric $\omega_1$ conformal to $\omega$ such that $S^{Ch}(\omega_1)= S$. Let $\tilde{\omega}=e^{\frac{2}{n}a}\omega_1$, then 
\begin{align}
S^{Ch}(\tilde{\omega})=e^{-\frac{2}{n}a}\left(-\Delta^{Ch}_{\omega_1} a+S^{Ch}(\omega_1)\right)=e^{-\frac{2}{n}a}S=\tilde{S}.
\end{align}
This completes the proof of the theorem under condition (1). 

Similarly, if 
\begin{align}
 0\leq \alpha\leq 2, \quad l=1-\frac{\alpha}{2}, \quad k= 1+\frac{\alpha}{2}, \quad b^2>4n^2\sqrt{C_3},
\end{align}
then we set 
\begin{align}\label{4.39-}
S=e^{\frac{2}{n}a}\tilde{S}, \quad \quad c^2=e^{\frac{2}{n}a}\tilde{c}^2,
\end{align}
where $a<\frac{n}{2}\log{\frac{b^2-4n^2\sqrt{C_3}}{\tilde{c}^2}}$ is a constant. From \eqref{4.33.} and \eqref{4.39-}, we have
\begin{align}
-c^2(1+r^k)\leq S<0, \quad c^2<b^2-4n^2\sqrt{C_3}.
\end{align}
Thus, by Proposition \ref{thm4.4} (3), there exists a complete Hermitian metric $\omega_2$ conformal to $\omega$ such that $S^{Ch}(\omega_2)= S$. Let $\tilde{\omega}=e^{\frac{2}{n}a}\omega_2$, then 
\begin{align}
S^{Ch}(\tilde{\omega})=e^{-\frac{2}{n}a}\left(-\Delta^{Ch}_{\omega_2} a+S^{Ch}(\omega_2)\right)=e^{-\frac{2}{n}a}S=\tilde{S}.
\end{align}
This finishes the proof of this theorem.
\qed

\proof[Proof of Theorem \ref{theorem1.1}]
Taking $k=0$ and $\tilde{S}\equiv-1$ in Theorem \ref{thm4.5}, the theorem follows.
\qed

\proof[Proof of Theorem \ref{theorem1.3}] Firstly, we reduce \eqref{1.8+} to the case where \eqref{1.8+} holds for all $M$. Let $\phi \in C^\infty_0(M)$ satisfy $\Delta^{Ch}_\omega \phi=\delta$ on $D$ and $\Delta^{Ch}_\omega \phi\geq -\frac{b^2}{2}$ on $M$, where $\delta>0$ is a constant. To construct such a $\phi$, let $D\subset D_1\subset\subset D_2$, where $D_1, D_2$ are relatively compact open subsets of $M$ with smooth boundaries. let $\phi_0\in C^\infty(\overline{D_2})$ be the solution of
\begin{align}
\begin{cases}
&\Delta^{Ch}_\omega \phi_0=\delta_1 \quad\text{in}\ D_2\\
& \phi_0|_{\partial D_2}=\delta_2
\end{cases}
\end{align}
where $\delta_1, \delta_2>0$ are constants. Extend $\phi_0|_{D_1}$ smoothly to a compactly supported function $\phi_1\in C^\infty_0(M)$. Then $\Delta^{Ch}_\omega \phi_1=\delta_1$ on $D$ and $\Delta^{Ch}_\omega \phi_1\geq-C$ on $M$ for some constant $C>0$. Defining $\phi=\frac{b^2}{2C}\phi_1\in  C^\infty_0(M)$, we get $\Delta^{Ch}_\omega \phi=\delta$ on $D$ and $\Delta^{Ch}_\omega \phi\geq -\frac{b^2}{2}$ on $M$, with $\delta:=\frac{b^2}{2C}\delta_1>0$. 

Let $\omega_1=e^{\frac{2}{n}\phi}\omega$, then
\begin{align}\label{4.60+}
S^{Ch}(\omega_1)=e^{-\frac{2}{n}\phi}\left(-\Delta^{Ch}_\omega\phi+S^{Ch}(\omega)\right)\leq -\epsilon \quad \text{in}\ M,
\end{align}
where we have used \eqref{1.7+} and \eqref{1.8+}, and $\epsilon=e^{-\frac{2}{n}\|\phi\|_{C^0(M)}}\min\{\delta, \frac{b^2}{2}\}$. If $\tilde{\omega}=e^{\frac{2}{n}u}\omega_1$ for some $u\in C^\infty(M)$, then $u$ satisfies
\begin{align}
-\Delta^{Ch}_{\omega_1} u+S^{Ch}(\omega_1)=S^{Ch}(\tilde{\omega})e^{\frac{2}{n}u}.
\end{align}
Now choose $u_-\equiv a$ where $a<\frac{n}{2}\log{\frac{\epsilon}{c^2}}$ is a constant. Then, by \eqref{4.60+} and \eqref{1.9}, we have
\begin{align}
-\Delta^{Ch}_{\omega_1} u_-+S^{Ch}(\omega_1)- Se^{\frac{2}{n}u_-}=S^{Ch}(\omega_1)-Se^{\frac{2}{n}a}\leq -\epsilon+c^2e^{\frac{2}{n}a}<0.
\end{align}
According to \eqref{1.9} and Proposition \ref{prop3.2}, there is a solution $\tilde{u}\in C^\infty(M)$ of 
\begin{align}
-\Delta^{Ch}_{\omega_1} u+S^{Ch}(\omega_1)={S}e^{\frac{2}{n}u}\quad \text{with}\ u\geq u_-=a.
\end{align}
Set $\tilde{\omega}=e^{\frac{2}{n}\tilde{u}}\omega_1$. Then $S^{Ch}(\tilde{\omega})={S}$. Since $\omega_1=e^{\frac{2}{n}\phi}\omega$, we have 
\begin{align}
\tilde{\omega}=e^{\frac{2}{n}(\tilde{u}+\phi)}\omega\geq e^{\frac{2}{n}(a-\|\phi\|_{C^0(M)})}\omega
\end{align}
and hence $\tilde{\omega}$ is complete.
\qed

\section*{Statements and Declarations}
\textbf{Data availability.} No data was used for the research described in the article.

\textbf{Conflict of interest.}  The author states that there is no conflict of interest. The author has no relevant financial or non-financial interests to disclose.

\bigskip

Weike Yu

School of Mathematical Sciences,

Ministry of Education Key Laboratory of NSLSCS,

Nanjing Normal University,

Nanjing, 210023, Jiangsu, P. R. China,

wkyu2018@outlook.com

\bigskip

\end{document}